\documentclass[11pt]{article}
\usepackage[T1]{fontenc}
\usepackage[utf8]{inputenc}
\usepackage{lmodern,amsmath,amssymb,amsthm,mathtools}
\usepackage[margin=1in]{geometry}
\usepackage[colorlinks=true,linkcolor=blue,citecolor=blue,urlcolor=blue]{hyperref}
\newtheorem{theorem}{Theorem}[section]
\newtheorem{lemma}[theorem]{Lemma}
\newtheorem{proposition}[theorem]{Proposition}
\newtheorem{corollary}[theorem]{Corollary}
\theoremstyle{remark}\newtheorem{remark}[theorem]{Remark}
\newcommand{\C}{\mathbb C}\newcommand{\R}{\mathbb R}\newcommand{\Z}{\mathbb Z}
\newcommand{\T}{\mathcal T}\newcommand{\Sspace}{\mathcal S}
\newcommand{\A}{\mathcal A}\newcommand{\Hset}{\mathcal H}
\newcommand{\M}{\mathcal M}
\DeclareMathOperator{\spanop}{span}\DeclareMathOperator{\rank}{rank}
\numberwithin{equation}{section}
\title{Obstructions to Hadamard Deformations of Fourier Matrices}
\author{Shuler Hopkins, University of Virginia\\Remus Nicoara, University of Tennessee Knoxville}\date{}
\begin{document}
\maketitle
\begin{abstract}
We study deformations of the Fourier spin model commuting square in the space of commuting squares, or
equivalently deformations of the matrix $\sqrt nF_n$, where $F_n$ denotes the Fourier matrix, by complex Hadamard matrices. If $n$ is neither a prime power nor
twice a prime, we show that there exist unit vectors in
the enveloping tangent space that cannot be realized as directions of
convergence of sequential deformations. 

Consequently, the derivative
of any differentiable Hadamard family through $\sqrt nF_n$ has rank
strictly smaller than the dimension of the enveloping tangent space at $\sqrt nF_n$ (i.e., the defect of $F_n$). This is particularly surprising because for every
$n$ there exist enough one-parameter analytic families to provide a basis of directions in the enveloping
tangent space. 

We obtain cubic obstructions when $n$ has at least three
distinct prime divisors, quartic obstructions when it has exactly two
and a repeated factor, and quintic obstructions when it is the product
of two distinct odd primes. These results extend our
order $30$ result from \cite{HN}.
\end{abstract}

\section{Introduction}
Commuting squares play an important role in the construction and
classification of subfactors, through the work of Jones \cite{Jo} and
Popa \cite{Po}. In the case of spin model commuting squares, their deformation theory is equivalent to the deformation theory of complex Hadamard matrices. Such deformations also have significant consequences to quantum information theory (see for instance \cite{TZ1}). In this paper
we study the local structure of these deformations at the Fourier matrix.

A complex Hadamard matrix is a matrix $H\in M_n(\C)$ with
$|H_{jk}|=1$ and $HH^*=nI_n$. We write $\Hset(n)$ for the set of
such matrices. Let $D_n\subset M_n(\C)$ be the diagonal algebra, let
\[
 F_n=\frac1{\sqrt n}(e^{2\pi i jk/n})_{0\le j,k<n},
 \qquad C_n=F_nD_nF_n^*,
\]
and denote by $\tau$ the normalized trace. For a unitary $U\in M_n(\C)$,
the following square of inclusions is a commuting square if and only if
$\sqrt n UF_n$ is a complex Hadamard matrix:
\[
 \begin{pmatrix}
 D_n&\subset&M_n(\C)\\
 \cup&&\cup\\
 \C I_n&\subset&UC_nU^*
 \end{pmatrix}.
\]
Equivalently,
\begin{equation}\label{eq:square}
 \tau(UpU^*q)=\tau(pq)\text{ for all }p\in C_n,\ q\in D_n.
\end{equation}
We use the inner product $\langle x,y\rangle=\tau(xy^*)$ and its associated
norm $\|x\|_2$. As usual, $[C_n,D_n]$ denotes the linear span of the
commutators $[p,q]=pq-qp$ with $p\in C_n,\ q\in D_n$.

The study of sequential deformations in \cite{N1,N3} gives necessary
conditions for a commuting square to be a limit of non-isomorphic
commuting squares. The first order conditions lead to a span condition
which implies isolation \cite{N1}; further restrictions appear at
second order \cite{N3}. For the Fourier commuting square, the first
order conditions have the following form.

A direction of convergence is a limit
\[
 a=\lim_{k\to\infty}\frac{U_k-I_n}{\|U_k-I_n\|_2},
\]
where $U_k\ne I_n$ are unitaries satisfying \eqref{eq:square} and
$U_k\to I_n$. Every such $a$ belongs to the real vector space
\[
 \A_n=\{a\in[C_n,D_n]^\perp:a^*=-a\}.
\]
Under the identification $a\mapsto\sqrt n aF_n$, this is the enveloping
tangent space of $\Hset(n)$ at $\sqrt n F_n$; see \cite{Ba,TZ2}. Its dimension, the undephased defect, is
\begin{equation}\label{eq:defect}
 d_n=\sum_{k=0}^{n-1}\gcd(n,k)
     =\prod_{p^e\parallel n}p^{e-1}\big((e+1)p-e\big),
 \text{ where }\gcd(n,0)=n.
\end{equation}
The dephased defect is $d_n'=d_n-2n+1$, corresponding to matrices whose
first row and column consist of ones; see \cite{NW1,TZ2}. These quantities
bound the number of independent tangent vectors of a differentiable family
through $\sqrt n F_n$. For prime-power orders, families attaining the
dephased bound are known \cite{TZ2}.

For every $n$, Nicoar\u a and White \cite{NW2} constructed $d_n$
one-parameter analytic families through $\sqrt nF_n$ whose tangent
vectors form a basis of the enveloping tangent space. It is therefore
natural to ask whether
these families can be combined into one differentiable family with $d_n$
independent parameters. The second order conditions do not rule this out:
they are satisfied by every element of $\A_n$ \cite{NW3}.

In \cite{HN} we showed that the answer is negative for
$n=30$. We obtained third order restrictions valid even for sequential
deformations, and constructed a vector in the enveloping tangent space that does not satisfy these conditions. The perturbative
computations of Barros e S\'a and Bengtsson \cite[Section~7]{BB}, for
$n\le100$ and also $n=210$, suggested third order obstructions for
orders with at least three distinct
prime divisors, fourth order obstructions for two-prime orders with
a repeated factor, and fifth order obstructions for products of two
distinct odd primes. We prove a uniform result covering these three cases.

\begin{theorem}\label{thm:main}
Let $n\ge2$ be neither a prime power nor twice a prime.
There exists a nonzero $a\in\A_n$ such that $a/\|a\|_2$ is not a
direction of convergence.
\end{theorem}

Consequently, the derivative of any differentiable Hadamard family
through $\sqrt nF_n$ has rank strictly less than $d_n$, and strictly less than
$d_n'$ in the dephased case (Corollary~\ref{cor:dimension}). The surprising
aspect is that a full basis of the enveloping tangent space is realized by analytic deformations (see  \cite{NW2}), while
some other directions in the enveloping tangent space cannot be realized even by convergent sequences.

We note that for prime-power orders it is known that there exist $d(F_n)$-parametric families that cover all directions in the enveloping tangent space (see \cite{TZ2}). Thus our theorem cannot be extended to that case. However, if $n=2q$ with
$q>3$ prime, we conjecture that the theorem is still true. When $q=3$ (so $n=6$) it is again not possible to extend the theorem: the dephased defect of $F_6$ is 4, and it was shown in \cite{Sz} that there exists a 4-parameter family of complex Hadamard matrices through $F_6$ that covers all directions in the enveloping tangent space at $F_6$.

We end this section by listing the main steps of the proof.

We start by deriving a cubic restriction directly from the commuting-square
equations; we do this in terms of sequential deformations, without stronger assumptions such as differentiability or analyticity.  When $n$ has
at least three distinct prime divisors, we write $n=rst$ with $r,s,t>1$
pairwise coprime and construct an explicit obstructed direction using
three tensor products of diagonal and shift unitaries. 

For $n=p^\alpha q^\beta$, with $p\ne q$ prime we show that the cubic conditions are automatically satisfied, then we look at stronger restrictions. We eliminate the variables orthogonal to the tangent space and obtain a quartic condition valid for arbitrary convergent
sequences. A nonzero coefficient gives an obstructed direction whenever
$\max(\alpha,\beta)\ge2$. 

Finally, for $n=pq$ with $p,q$ distinct odd
primes, we prove that the quartic conditions always hold then
show that a fifth degree condition fails.

\medskip
\noindent\textbf{AI Statement.}
ChatGPT Astra was used in extending some of the arguments,
performing symbolic computations, preparing some parts of the manuscript, and formalizing the manuscript in
Lean~4. The authors take full responsibility for the
mathematical content and the final version of the paper.

\section{A third order condition on sequential deformations}\label{sec:third}
In this section we prove third order restrictions on any direction of convergence of a sequence of Hadamard matrices converging to $F_n$, and different from $F_n$.

 Let $n=n_1\cdots n_k$, where the $n_\ell>1$ are pairwise coprime.
The Chinese remainder theorem identifies the Fourier commuting square
with the tensor product square given by
\[
 F=\bigotimes_{\ell=1}^k F_{n_\ell},\qquad
 D=\bigotimes_{\ell=1}^k D_{n_\ell},\qquad
 C=\bigotimes_{\ell=1}^k C_{n_\ell}.
\]
The identification is by a permutation of the standard basis. Put
$\Sspace=[C,D]$, $\T=\Sspace^\perp$, and
$\A=\{a\in\T:a^*=-a\}$.
For $m>1$ let
\[
 P_m=\sum_{j\in\Z_m}e_{j,j-1},\qquad
 Q_m=\sum_{j\in\Z_m}e^{2\pi i j/m}e_{j,j}.
\]
For $x,y\in G=\prod_{\ell=1}^k\Z_{n_\ell}$, set
\[
 Q_x=\bigotimes_{\ell=1}^k Q_{n_\ell}^{x_\ell},\qquad
 P_y=\bigotimes_{\ell=1}^k P_{n_\ell}^{y_\ell},\qquad
 \chi(x,y)=\exp\left(2\pi i\sum_{\ell=1}^k\frac{x_\ell y_\ell}{n_\ell}\right).
\]
Write $E_{x,y}=Q_xP_y$. These matrices form an orthonormal basis of
$M_n(\C)$, with
\begin{equation}\label{eq:weyl}
 E_{x,y}E_{z,h}=\chi(z,y)^{-1}E_{x+z,y+h},\qquad
 \tau(E_{x,y})=\delta_{x,0}\delta_{y,0}.
\end{equation}
Since $[P_y,Q_x]=(\chi(x,y)^{-1}-1)E_{x,y}$, we have
\begin{equation}\label{eq:spaces}
 \T=\spanop\{E_{x,y}:\chi(x,y)=1\},\qquad
 \Sspace=\spanop\{E_{x,y}:\chi(x,y)\ne1\}.
\end{equation}
The matrices $P_y$ and $Q_x$ form bases of $C$ and $D$, respectively.
We call $(P_i,Q_j)$ a basis pair; it commutes exactly when $\chi(i,j)=1$.
Both spaces in \eqref{eq:spaces} are closed under adjoints, so $\A$ is a
real form of $\T$: every element of $\T$ is uniquely $a+ib$ with
$a,b\in\A$.
Coprimality implies that
$\chi(x,y)=1$ exactly when $x_\ell y_\ell=0$ in each cyclic coordinate.
For such indices, $E_{x,y}^*=E_{-x,-y}$.

\begin{lemma}\label{lem:vanishing}
If $X,Y\in\T$ and $\chi(i,j)=1$, then $\tau([P_i,X]YQ_j)=0$.
\end{lemma}
\begin{proof}
Recall that in $\Z_m$ the relations
$xy=zh=(x+z)(y+h)=0$ imply $xh=zy=0$ \cite{HN}.
Indeed, work modulo $p^e\mid m$, using valuations truncated at $e$.
If $v_p(x)\le v_p(z)$,
then $xy=0$ implies $zy=0$; otherwise $zh=0$ implies $xh=0$.
Since $xh+zy=0$, both vanish.

It suffices to take $X=E_{x,y}$ and $Y=E_{z,h}$ from the basis of $\T$.
By \eqref{eq:weyl}, the trace is zero unless $x+z=-j$ and $y+h=-i$.
In that case it equals $\chi(x,h)-1$, which is zero by the preceding
observation applied in each coordinate.
\end{proof}

For $a=\sum_{\chi(x,y)=1}\alpha_{x,y}E_{x,y}\in\T$, define
\begin{equation}\label{eq:B}
 B(a)=\sum_{\substack{\chi(x,y)=\chi(z,h)=1\\\chi(x+z,y+h)\ne1}}
 \alpha_{x,y}\alpha_{z,h}
 \frac{1-\chi(x,h)}{1-\chi(x+z,y+h)}E_{x+z,y+h}.
\end{equation}
This is the unique element of $\Sspace$ satisfying
\begin{equation}\label{eq:second}
 \tau(B(a)[p,q])=-\tau([p,a]a q)\qquad(p\in C,\ q\in D).
\end{equation}
Indeed, for $i=-(y+h)$ and $j=-(x+z)$, the trace pairing of
$E_{x+z,y+h}$ with $[P_i,Q_j]$ is $1-\chi(x+z,y+h)$, whereas the
ordered pair $(x,y),(z,h)$ contributes $1-\chi(x,h)$ on the right.
The remaining equations follow from Lemma~\ref{lem:vanishing}.
Uniqueness follows from \eqref{eq:spaces}.

\begin{proposition}\label{prop:cubic}
For $a\in\T$ and $\chi(i,j)=1$, put
\begin{equation}\label{eq:cubic}
 \Phi_{i,j}(a)=\tau([P_i,a](a^2-B(a))Q_j)
                 -\tau([P_i,B(a)]aQ_j).
\end{equation}
If $a$ is a direction of convergence, then $\Phi_{i,j}(a)=0$.
\end{proposition}
\begin{proof}
Let $V_k=U_k-I$, $\varepsilon_k=\|V_k\|_2$, and suppose
$V_k/\varepsilon_k\to a$. Unitarity and \eqref{eq:square} give
\begin{equation}\label{eq:exact}
 V_k+V_k^*=-V_kV_k^*,\qquad
 \tau(V_k[p,q])=\tau([p,V_k]V_k^*q).
\end{equation}
In particular $a\in\A$. Write $V_k=X_k+Y_k$, with
$X_k\in\T$ and $Y_k\in\Sspace$. The second identity gives
$Y_k=O(\varepsilon_k^2)$, since the commutators span $\Sspace$ and
its trace pairing is nondegenerate. Dividing that identity by
$\varepsilon_k^2$ and using \eqref{eq:second}, we obtain
\begin{equation}\label{eq:limits}
 X_k/\varepsilon_k\longrightarrow a,\qquad
 Y_k/\varepsilon_k^2\longrightarrow B(a).
\end{equation}
Projection of the first identity in \eqref{eq:exact} onto $\Sspace$ gives
\[
 B(a)+B(a)^*=\operatorname{proj}_{\Sspace}(a^2).
\]
Now let $p=P_i$, $q=Q_j$ with $[p,q]=0$. Then
$\tau([p,V_k]V_k^*q)=0$, while
$\tau([p,X_k]X_k^*q)=0$ by Lemma~\ref{lem:vanishing}.
After division by $\varepsilon_k^3$, the remaining terms converge to
\[
 \tau([p,a]B(a)^*q)-\tau([p,B(a)]a q)=0;
\]
the term quadratic in $Y_k$ tends to zero. Since
$B(a)^*-(a^2-B(a))\in\T$, another application of
Lemma~\ref{lem:vanishing} proves \eqref{eq:cubic}.
\end{proof}

\section{The case of three or more distinct prime divisors}\label{sec:three}

In this section we use the third order restrictions obtained in the previous section to prove the following:

\begin{proposition}\label{prop:three}
If $n$ has at least three distinct prime divisors, there exists a nonzero
$a\in\A_n$ such that $a/\|a\|_2$ is not a direction of convergence.
\end{proposition}
\begin{proof}
Group the full prime-power factors of $n$ to write $n=rst$, where
$r,s,t>1$ are pairwise coprime, and relabel so that $s>2$.
Use the tensor coordinates of Section~\ref{sec:third}, put
$u=e^{2\pi i/r}$, $v=e^{2\pi i/s}$, $w=e^{2\pi i/t}$, and set
\begin{equation}\label{eq:XYZ}
 \begin{aligned}
 X&=Q_r\otimes P_s\otimes P_t,\\
 Y&=Q_r^{-1}\otimes Q_s\otimes P_t^{-1},\\
 Z&=P_r\otimes Q_s^{-1}\otimes Q_t,
 \end{aligned}
 \qquad a=X+Y+Z-X^*-Y^*-Z^*.
\end{equation}
Each term belongs to $\T$. The three pairs of opposite Weyl indices
have different patterns of zero coordinates, and $X,Y,Z$ each have order
$rst>2$. Thus the six terms are distinct and $0\ne a\in\A$.
Let
\[
 P=P_r^{-1}\otimes P_s^{-1}\otimes I_t,\qquad
 Q=I_r\otimes I_s\otimes Q_t^{-1}.
\]
These commute. We will show that the corresponding expression
\eqref{eq:cubic} satisfies
\begin{equation}\label{eq:value}
 \Phi(a)=-\frac{(u-1)^2(v-1)(w-1)^2}{(uw-1)(uw-v)}\ne0.
\end{equation}
Since $B$ is quadratic, $\Phi$ is homogeneous of degree three.
Proposition~\ref{prop:cubic} will then exclude $a/\|a\|_2$.

The indices of $X,Y,Z$ are, respectively,
\[
 \xi_X=((1,0,0),(0,1,1)),\quad
 \xi_Y=((-1,1,0),(0,0,-1)),\quad
 \xi_Z=((0,-1,1),(1,0,0)).
\]
Their sum is $((0,0,1),(1,1,0))$, the negative of the index of $PQ$.
We claim that the only triples from $X,Y,Z,X^*,Y^*,Z^*$ with this sum
are the six permutations of $X,Y,Z$. Let $k_X,k_Y,k_Z$ be the signed
numbers of occurrences in such a triple. The $P_r$ and $Q_t$ coordinates
give $k_Z\equiv1\pmod r$ and $k_Z\equiv1\pmod t$.
Since $|k_Z|\le3$ and $rt\ge6$, we have $k_Z=1$.
The $P_s$ and $Q_s$ coordinates give $k_X\equiv k_Y\equiv1\pmod s$.
Now $|k_X|+|k_Y|+|k_Z|\le3$ and $s>2$ imply $k_X=k_Y=1$,
proving the claim.

Write $E_\xi=E_{x,y}$ for $\xi=(x,y)$.
For $L\in\{X,Y,Z\}$, let $\eta_L$ be the sum of the other two indices,
and put
\[
 d_L=\tau([P,L]E_{\eta_L}Q),\qquad
 e_L=\tau([P,E_{\eta_L}]LQ).
\]
Let $c_L$ and $b_L$ be the coefficients of $E_{\eta_L}$ in $a^2$ and
$B(a)$, respectively. The preceding claim shows that only the two orders
of the complementary pair contribute to these coefficients. From
\eqref{eq:weyl} and \eqref{eq:B} we obtain
\begin{equation}\label{eq:table}
\renewcommand{\arraystretch}{1.8}
\begin{array}{c|cccc}
 L&d_L&e_L&c_L&b_L\\\hline
 X&\dfrac{u-1}{uw}&-\dfrac{u-1}{u}&u+w&
       \dfrac{2uw-u-w}{uw-1}\\
 Y&-\dfrac{w(u-v)}v&\dfrac{u-v}v&\dfrac vw+\dfrac1u&
       \dfrac{uv+w-2v}{uw-v}\\
 Z&1-v&v-1&1+\dfrac1v&1
\end{array}.
\end{equation}
For example, the ordered pairs $(Y,Z)$ and $(Z,Y)$ contribute
$w(u-1)/(uw-1)$ and $u(w-1)/(uw-1)$ to $b_X$.
The trace entries follow from
\[
 \tau([P_i,E_{x,y}]E_{z,h}Q_j)
 =\frac{\chi(i,x)^{-1}-1}{\chi(i+y,z)}
 \quad\text{when }(x+z,y+h)=(-j,-i).
\]
Substituting \eqref{eq:table} into \eqref{eq:cubic} and collecting terms gives
\[
 \Phi(a)=\sum_{L=X,Y,Z}\bigl(d_Lc_L-(d_L+e_L)b_L\bigr)
       =-\frac{(u-1)^2(v-1)(w-1)^2}{(uw-1)(uw-v)}.
\]
The numerator is nonzero. Also, $uw$ has order $rt$, while $v$ has
order $s$, so coprimality gives $uw\ne1$ and $uw\ne v$.
This proves \eqref{eq:value}, and hence the proposition.
\end{proof}

\section{A fourth order condition on sequential deformations}\label{sec:fourth}

In this section we show that when $n$ has just two distinct prime divisors the third order conditions do not yield a contradiction, then we proceed with developing fourth order conditions. Throughout the section, let $n=p^\alpha q^\beta$, where $p,q$ are
distinct primes and $\alpha,\beta\ge1$. We use the notation of
Section~\ref{sec:third} with $m=p^\alpha$, $r=q^\beta$, and
$G=\Z_m\times\Z_r$.

\subsection{Reduction to the tangent space}

We use the equations corresponding to non-commuting basis pairs
to express the component orthogonal to $\T$ as a holomorphic function
of the component in $\T$. Substituting this expression into the
remaining equations gives conditions involving only variables in
$\T$, from which we obtain the higher order obstructions.

For $V$ sufficiently close to zero, consider the holomorphic equations
\begin{equation}\label{eq:holomorphic}
 \tau((I+V)p_0(I+V)^{-1}q_0)=\tau(p_0q_0),\qquad p_0\in C,\ q_0\in D.
\end{equation}
Every unitary Hadamard deformation satisfies these equations. Expansion
of the inverse gives the equivalent identities
\begin{equation}\label{eq:inverse}
 \tau(V[p_0,q_0])=\sum_{j\ge1}(-1)^j\tau([p_0,V]V^j q_0).
\end{equation}
For $p_0=P_{-y}$ and $q_0=Q_{-x}$, the coefficient of $E_{x,y}$ on the
left is $1-\chi(x,y)$. Thus the equations indexed by
$\chi(x,y)\ne1$ have an invertible derivative on $\Sspace$.

The holomorphic implicit function theorem gives a unique holomorphic
map $Y$ from a neighborhood of zero in $\T$ to $\Sspace$, with
$Y(0)=DY(0)=0$, such that $V=X+Y(X)$ satisfies the equations indexed
by $\chi(x,y)\ne1$. Conversely, every sufficiently small solution of
these equations has this form. Write
\begin{equation}\label{eq:Y}
 Y(X)=B(X)+K(X)+O(\|X\|_2^4),
\end{equation}
where $B$ and $K$ are homogeneous of degrees two and three.
The quadratic term is exactly the map $B$ in \eqref{eq:B}. The cubic
term is determined, for the non-commuting basis pairs, by
\begin{equation}\label{eq:K}
 \tau(K(X)[p_0,q_0])=
 \tau([p_0,X](X^2-B(X))q_0)-\tau([p_0,B(X)]Xq_0).
\end{equation}
For a commuting basis pair, the quadratic term on the right of
\eqref{eq:inverse} vanishes identically by Lemma~\ref{lem:vanishing}.
We next show that the cubic term vanishes as well.

\subsection{Absence of third order obstructions}

We show that if $n$ has only two distinct prime divisors then the third order conditions are always satisfied. This justifies looking for higher order restrictions.

\begin{lemma}\label{lem:star}
Let $\ell$ be a prime power and suppose that
$x_i y_i=0$ for $1\le i\le3$, and
$(x_1+x_2+x_3)(y_1+y_2+y_3)=0$ in $\Z_\ell$.
The off-diagonal entries of the matrix
$(e^{2\pi i x_i y_j/\ell})_{i,j=1}^3$ are all one, except possibly
for two entries in one row or one column, whose values are $z,z^{-1}$.
\end{lemma}
\begin{proof}
Write $\ell=p^e$ and use valuations truncated at $e$, with $v_p(0)=e$. Let $a$ and $b$
be the minimum valuations of the $x_i$ and $y_i$. If $a+b\ge e$,
all cross products vanish. Otherwise the two sets of minimizing indices
are disjoint, since $x_i y_i=0$. They cannot both be singletons:
then the valuations of the two sums would be $a$ and $b$.
Thus one set has two indices and the other has one. If two $x_i$
have valuation $a$, their corresponding $y_i$ have valuation at least
$e-a$, so only the two cross products in the remaining column can be
nonzero. Their sum is zero by the hypothesis on the two sums.
The case of two minimizing $y_i$ is the same with rows and columns
interchanged.
\end{proof}

\begin{proposition}\label{prop:cubiczero}
For $X\in\T$ and a commuting basis pair $p_0,q_0$,
\begin{equation}\label{eq:cubiczero}
 \tau([p_0,X](X^2-B(X))q_0)-\tau([p_0,B(X)]Xq_0)=0.
\end{equation}
\end{proposition}
\begin{proof}
Let $W_i=E_{x_i,y_i}\in\T$, $1\le i\le3$. It suffices to show that
the coefficient of $z_1z_2z_3$ in the left-hand side of
\eqref{eq:cubiczero}, evaluated at $X=\sum_i z_iW_i$, is zero.
The matrices $W_i$ need not be distinct. This coefficient vanishes
unless the sum of their indices is the negative of the index of
$p_0q_0$. In the remaining case set
\[
 t_{ij}=\chi(x_i,y_j),\quad t_{ii}=1,\quad
 \rho_i=\prod_{j=1}^3t_{ij},\quad \sigma_{jk}=t_{jk}t_{kj}.
\]
Here $\prod_i\rho_i=1$. Define
\[
 b_{jk}=\begin{cases}
 (2-t_{jk}-t_{kj})/(1-\sigma_{jk}),&\sigma_{jk}\ne1,\\
 0,&\sigma_{jk}=1.
 \end{cases}
\]
By \eqref{eq:B} and \eqref{eq:weyl}, this coefficient is
\begin{equation}\label{eq:cubiccoefficient}
 \sum_{i=1}^3
 \left[(\rho_i-1)\sigma_{jk}(t_{jk}^{-1}+t_{kj}^{-1})
 -\big((\rho_i-1)\sigma_{jk}+\rho_i^{-1}-1\big)b_{jk}\right].
\end{equation}
Here $\{j,k\}=\{1,2,3\}\setminus\{i\}$.

Apply Lemma~\ref{lem:star} in each prime-power coordinate. Each prime
contributes either no nontrivial entry, a row pair $R_i(z)$, or a
column pair $C_i(z)$. In $R_1(z)$ the entries $t_{12},t_{13}$ are
$z,z^{-1}$; define the other pairs by using increasing order for the
two remaining indices. A single row or column pair gives zero in
\eqref{eq:cubiccoefficient}. Up to a simultaneous permutation of the
three indices, inversion of parameters, and interchange of the two
prime factors, all cases with two nontrivial pairs appear below.
Direct substitution gives zero in every row.
\begin{center}
\renewcommand{\arraystretch}{1.35}
\begin{tabular}{c|cccccc}
 pairs&$t_{12}$&$t_{13}$&$t_{21}$&$t_{23}$&$t_{31}$&$t_{32}$\\\hline
 $R_1(u),R_1(v)$&$uv$&$(uv)^{-1}$&$1$&$1$&$1$&$1$\\
 $R_1(u),R_2(v)$&$u$&$u^{-1}$&$v$&$v^{-1}$&$1$&$1$\\
 $R_1(u),C_1(v)$&$u$&$u^{-1}$&$v$&$1$&$v^{-1}$&$1$\\
 $R_1(u),C_2(v)$&$uv$&$u^{-1}$&$1$&$1$&$1$&$v^{-1}$\\
 $C_1(u),C_1(v)$&$1$&$1$&$uv$&$1$&$(uv)^{-1}$&$1$\\
 $C_1(u),C_2(v)$&$v$&$1$&$u$&$1$&$u^{-1}$&$v^{-1}$
\end{tabular}
\end{center}
There are no additional denominator degeneracies: $u$ and $v$ are
nontrivial roots of coprime prime-power orders, so a product
$u^\varepsilon v^\delta$, with $\varepsilon,\delta\in\{-1,0,1\}$,
is one only if both exponents vanish. The cases $u=1$ or $v=1$
were already included by omitting the corresponding prime contribution. This proves the proposition.
\end{proof}

\subsection{The fourth order condition}

In this section we develop fourth order obstructions for possible directions of convergence of sequential families of complex Hadamard matrices.

For a commuting basis pair $p_0,q_0$, write the remaining equation as
\begin{equation}\label{eq:R}
 \mathcal R_{p_0,q_0}(X)=\tau(p_0q_0)-
 \tau((I+X+Y(X))p_0(I+X+Y(X))^{-1}q_0).
\end{equation}
Its terms of degrees two and three vanish by Lemma~\ref{lem:vanishing}
and Proposition~\ref{prop:cubiczero}. Expansion of \eqref{eq:inverse}
therefore gives
\begin{equation}\label{eq:Rquartic}
 \mathcal R_{p_0,q_0}(X)=\Psi_{p_0,q_0}(X)+O(\|X\|_2^5),
\end{equation}
where the homogeneous quartic polynomial is
\begin{align}\label{eq:Psi}
 \Psi_{p_0,q_0}(X)={}&-\tau([p_0,X]K(X)q_0)-\tau([p_0,K(X)]Xq_0)
                 -\tau([p_0,B(X)]B(X)q_0)\notag\\
 &+\tau([p_0,B(X)]X^2q_0)
  +\tau([p_0,X](XB(X)+B(X)X)q_0)-\tau([p_0,X]X^3q_0).
\end{align}
\begin{proposition}\label{prop:necessary}
If $n$ has exactly two distinct prime divisors, every direction of
convergence $a$ satisfies $\Psi_{p_0,q_0}(a)=0$ for every commuting basis pair.
\end{proposition}
\begin{proof}
Write $V_k=U_k-I=X_k+Y_k$ with $X_k\in\T$ and $Y_k\in\Sspace$.
The implicit function theorem gives $Y_k=Y(X_k)$ for all sufficiently
large $k$. For a commuting pair, \eqref{eq:Rquartic} gives
\[
 0=\Psi_{p_0,q_0}(X_k)+O(\|X_k\|_2^5).
\]
If $\varepsilon_k=\|V_k\|_2$ and $V_k/\varepsilon_k\to a$, then
$X_k/\varepsilon_k\to a$. Since $\|X_k\|_2\le\varepsilon_k$, division by $\varepsilon_k^4$ gives
$\Psi_{p_0,q_0}(a)=0$. No differentiability of the sequence is required.
\end{proof}

\section{The case of two prime divisors with a repeated factor}\label{sec:two}
We retain the notation of Section~\ref{sec:fourth}. The following
calculation will provide a nonzero quartic obstruction for $n=p^\alpha q^\beta$, where $p\ne q$ are primes,
$\alpha,\beta\ge1$, and $\max(\alpha,\beta)\ge2$.

\begin{lemma}\label{lem:coefficient}
Let $W_i=E_{x_i,y_i}\in\T$, $1\le i\le4$, and suppose their cross
characters are
\begin{equation}\label{eq:matrix}
 (\chi(x_i,y_j))_{i,j=1}^4=
 \begin{pmatrix}
 1&1&v&v^{-1}\\
 1&1&v&v^{-1}\\
 \lambda&\mu&1&1\\
 \nu&(\lambda\mu\nu)^{-1}&1&1
 \end{pmatrix}.
\end{equation}
For $I\subset\{1,2,3,4\}$, write $x_I=\sum_{i\in I}x_i$ and
$y_I=\sum_{i\in I}y_i$. Assume that $\chi(x_I,y_I)\ne1$ when
$|I|=3$, or when $|I|=2$ and $I\ne\{1,2\},\{3,4\}$.
Put $p_0=P_{-\sum y_i}$ and
$q_0=Q_{-\sum x_i}$. These commute, and the coefficient of
$z_1z_2z_3z_4$ in $\Psi_{p_0,q_0}(\sum_i z_iW_i)$ is
\begin{equation}\label{eq:factor}
 \frac{(v-1)^4(\lambda\mu-1)^2(\lambda\nu-1)^2}
 {\lambda(v-\nu)(\lambda v-1)(\mu v-1)(\lambda\mu\nu v-1)}.
\end{equation}
The matrices $W_i$ need not be distinct.
\end{lemma}
\begin{proof}
Put $X=\sum_i z_iW_i$ and $\theta_I=\chi(x_I,y_I)$. We compute the
coefficient using \eqref{eq:weyl}. For disjoint subsets define
\[
 m_{I,J}=\chi(x_J,y_I)^{-1},\qquad
 d_{I,J}=(\chi(x_I,y_I+y_J)-1)\theta_J.
\]
These are, respectively, the multiplication coefficient and the trace
coefficient
\[
 \tau([P_{-y_I-y_J},E_{x_I,y_I}]E_{x_J,y_J}Q_{-x_I-x_J})=d_{I,J}.
\]
Let $a_I,b_I$ be the coefficients of $\prod_{i\in I}z_i$ in $X^2,B(X)$,
relative to $E_{x_I,y_I}$, for $|I|=2$. Then
\[
 a_{ij}=t_{ij}^{-1}+t_{ji}^{-1},\qquad
 b_{12}=b_{34}=0,\qquad
 b_{ij}=\frac{2-t_{ij}-t_{ji}}{1-t_{ij}t_{ji}}\quad(ij\ne12,34),
\]
where $t_{ij}=\chi(x_i,y_j)$.
For $|I|=3$, the coefficient $k_I$ in $K(X)$ is
\begin{equation}\label{eq:kI}
 k_I=\frac{\sum_{i\in I}
   \{d_{i,J}(a_J-b_J)-d_{J,i}b_J\}}{1-\theta_I},
 \qquad J=I\setminus\{i\}.
\end{equation}
Singleton braces in subscripts are suppressed. Only $k_{123}$ and
$k_{124}$ will be needed: for $i=1,2$ and $J=\{1,2,3,4\}\setminus\{i\}$,
the first two row products of \eqref{eq:matrix} give $d_{i,J}=d_{J,i}=0$.
Set $\delta=(\lambda\mu\nu)^{-1}$ and
\[
 L(a,b,w)=\frac{(w-1)^2\{w[(a+b)(ab+1)-4ab]-2(a-1)(b-1)\}}
 {(aw-1)(bw-1)(abw^2-1)}.
\]
Substitution in \eqref{eq:kI} gives
\[
 k_{123}=L(\lambda,\mu,v),\qquad k_{124}=L(\nu,\delta,v^{-1}).
\]
To finish, put, for $|I|=3$,
\[
 r_I=\sum_{i\in I}(m_{i,J}+m_{J,i})b_J,\qquad
 s_I=\sum_{i\in I}m_{i,J}a_J,\qquad J=I\setminus\{i\}.
\]
These are the coefficients in $XB+BX$ and $X^3$.
With $S=\{1,2,3,4\}$, expansion of \eqref{eq:Psi} therefore gives
\begin{align}\label{eq:finite}
 &\sum_{i\in S}\left[-(d_{i,J}+d_{J,i})k_J+d_{i,J}(r_J-s_J)\right]
       +\sum_{\substack{I\subset S\\|I|=2}}d_{I,S\setminus I}
                     b_I(a_{S\setminus I}-b_{S\setminus I}),
\end{align}
where $J=S\setminus\{i\}$ in the first sum.
The first two terms of the first sum vanish, as noted above. Insert the displayed
values of $a_I,b_I,k_I,r_I,s_I$ into \eqref{eq:finite}, replace
$\delta$ by $(\lambda\mu\nu)^{-1}$, and put the fractions over a common denominator.
The factors $\lambda\mu v^2-1$ cancel; the result is \eqref{eq:factor}.
This is an identity of rational functions, so it holds at every
specialization satisfying the stated denominator conditions.
\end{proof}

\begin{proposition}\label{prop:two}
Let $n=p^\alpha q^\beta$, where $p\ne q$ are primes,
$\alpha,\beta\ge1$, and $\max(\alpha,\beta)\ge2$.
There exists a nonzero $a\in\A_n$ such that $a/\|a\|_2$ is not a
direction of convergence.
\end{proposition}
\begin{proof}
Relabel the primes so that $\alpha\ge2$. Put $m=p^\alpha$,
$r=q^\beta$, $h=p^{\alpha-1}$, and choose
\begin{equation}\label{eq:witness}
 \begin{aligned}
 W_1&=P_m\otimes Q_r,&W_2&=P_m^{h-1}\otimes Q_r,\\
 W_3&=Q_m\otimes P_r,&W_4&=Q_m^{p-1}\otimes P_r^{-1}.
 \end{aligned}
\end{equation}
All four matrices belong to $\T$. Let $\zeta=e^{2\pi i/m}$ and
$v=e^{2\pi i/r}$. Their cross characters have the form
\eqref{eq:matrix}, with
\[
 \lambda=\zeta,\qquad \mu=\zeta^{h-1},\qquad \nu=\zeta^{p-1}.
\]
Indeed, $\zeta^{(p-1)(h-1)}=(\lambda\mu\nu)^{-1}$ because $ph=m$.
The total Weyl index is $((p,2),(h,0))$, so we use the commuting pair
\[
 p_0=P_m^{-h}\otimes I_r,\qquad
 q_0=Q_m^{-p}\otimes Q_r^{-2}.
\]
The four cross-pair self-characters are $\lambda v,\mu v,\nu/v,(\lambda\mu\nu v)^{-1}$.
They are nontrivial: $\lambda,\mu,\nu,\lambda\mu\nu$ are roots of $p$-power order
different from one, whereas $v$ has $q$-power order.
The four triple self-characters are
\[
 \lambda\mu v^2,\quad(\lambda\mu v^2)^{-1},\quad\lambda\nu,\quad(\lambda\nu)^{-1}.
\]
They too are nontrivial, since $\lambda\mu=\zeta^h$ has order $p$ and
$\lambda\nu=\zeta^p$ has order $p^{\alpha-1}$. Thus Lemma~\ref{lem:coefficient}
applies, and its numerator and denominator are nonzero.
It follows that $\Psi_{p_0,q_0}$ is a nonzero polynomial on $\T$.

Since $\A$ is a real form of $\T$, this polynomial cannot vanish
identically on $\A$. Explicitly, choose a real basis of $\A$; it is
also a complex basis of $\T$, and a complex polynomial vanishing at
every real coordinate tuple is identically zero. Hence some
$0\ne a\in\A$ satisfies $\Psi_{p_0,q_0}(a)\ne0$.
By homogeneity and Proposition~\ref{prop:necessary}, $a/\|a\|_2$
cannot be a direction of convergence. The Chinese remainder
identification carries this conclusion back to $F_n$.

\end{proof}

\begin{remark}
For $n=12$, take $m=4,r=3$, so $\lambda=\mu=\nu=i$ and
$v=e^{2\pi i/3}$. The coefficient \eqref{eq:factor} is $-144$.
The same factorization applies to every pair of distinct primes and
all positive exponents with at least one greater than one.
\end{remark}

\section{Products of two distinct odd primes}\label{sec:odd}
In this section we deal with the case $n=pq$, where $p,q$ are distinct odd primes.
We retain the notations of Section~\ref{sec:fourth}, now with
$G=\Z_p\times\Z_q$. We first show that the quartic conditions vanish
identically, and then obtain a nonzero fifth order condition.

\subsection{Computing higher order conditions}

We derive recursive formulas for the coefficients of the reduced
equations from the previous sections. These formulas will be used to show that the fourth order
conditions vanish identically and to obtain a fifth order obstruction.

Take $W_i=E_{x_i,y_i}\in\T$, $1\le i\le k$, and put
\[
 t_{ij}=\chi(x_i,y_j),\qquad
 \theta_I=\prod_{i,j\in I}t_{ij},\qquad
 m_{I,J}=\prod_{i\in I,\,j\in J}t_{ji}^{-1}.
\]
Here $I,J$ are subsets of $\{1,\ldots,k\}$. Write
$x_I=\sum_{i\in I}x_i$, $y_I=\sum_{i\in I}y_i$,
$E_I=E_{x_I,y_I}$, and $z_I=\prod_{i\in I}z_i$.
For $X=\sum_i z_iW_i$, let $b_I E_I$ and $h_I E_I$ be the coefficients
of $z_I$ in $V=X+Y(X)$ and $(I+V)^{-1}$, respectively. Only square-free
monomials are needed, even when some $W_i$ coincide.
Set $b_{\{i\}}=1$, $h_{\{i\}}=-1$. For $|I|\ge2$, define
\begin{align}
 H_I&=-\sum_{J\sqcup K=I}b_Jh_Km_{J,K},\label{eq:H}\\
 r_I&=-\theta_I H_I-\sum_{J\sqcup K=I}b_Jh_K\theta_K,\label{eq:r}\\
 b_I&=\begin{cases}r_I/(1-\theta_I),&\theta_I\ne1,\\0,&\theta_I=1,
 \end{cases}
 \qquad h_I=H_I-b_I.\label{eq:b}
\end{align}
All sums here are over ordered decompositions into nonempty disjoint
subsets. The zero branch in \eqref{eq:b} follows from $Y(X)\in\Sspace$:
for $|I|\ge2$, the coefficient $b_IE_I$ is zero when $E_I\in\T$.
If $\theta_I=1$, then $r_I$ is the coefficient of $z_I$ in
$\mathcal R_{P_{-y_I},Q_{-x_I}}(X)$.
Indeed, \eqref{eq:H} follows by multiplying $I+V$ by its inverse, and
\eqref{eq:weyl} gives
\[
 \tau(E_JP_{-y_I}E_KQ_{-x_I})=\theta_K
 \quad(J\sqcup K=I).
\]
The coefficient in \eqref{eq:holomorphic} is
$(1-\theta_I)b_I-r_I$, proving the recurrence.

\subsection{Absence of quartic obstructions}
Since $p$ and $q$ are prime and distinct, $\chi(x,y)=1$ if and only if
$x_1y_1=0$ in $\Z_p$ and $x_2y_2=0$ in $\Z_q$. Consequently,
\begin{equation}\label{eq:mixed}
 \T=(C+D)\oplus\M,
\end{equation}
where $\M$ is spanned by
\[
 \mathsf A_{a,b}=P_p^a\otimes Q_q^b,\qquad
 \mathsf B_{c,d}=Q_p^c\otimes P_q^d,
 \quad a,c\ne0\pmod p,\ b,d\ne0\pmod q.
\]

\begin{lemma}\label{lem:gauge}
If $\Psi_{P,Q}$ vanishes on $\M$ for every commuting basis pair $P,Q$,
then it vanishes on $\T$ for every such pair.
\end{lemma}
\begin{proof}
For $Z\in\M$, the curve $U(t)=I+tZ+Y(tZ)$ satisfies all equations
\eqref{eq:holomorphic} modulo $t^5$, by the hypothesis and \eqref{eq:Rquartic}.
For $d\in D$ and $c\in C$, the curve
$\widetilde U(t)=e^{td}U(t)e^{tc}$ has exactly the same residuals,
because $c$ commutes with every $P\in C$ and $d$ with every $Q\in D$.
Its first derivative is $Z+d+c$.

Write $\widetilde U(t)-I=X(t)+N(t)$ according to $\T\oplus\Sspace$.
The noncommuting equations are satisfied exactly, so uniqueness in the
implicit function theorem gives $N(t)=Y(X(t))$. Thus \eqref{eq:Rquartic} gives
$\Psi_{P,Q}(X(t))=O(t^5)$, with $X(t)=t(Z+d+c)+O(t^2)$.
Dividing by $t^4$ and using \eqref{eq:mixed} proves the assertion.
\end{proof}

\begin{proposition}\label{prop:fourzero}
For every commuting basis pair $P,Q$, the polynomial $\Psi_{P,Q}$
vanishes identically on $\T$.
\end{proposition}
\begin{proof}
By Lemma~\ref{lem:gauge}, it suffices to work on $\M$. Polarize
$\Psi_{P,Q}$ and choose each $W_i$, $1\le i\le4$, of type
$\mathsf A$ or $\mathsf B$.
A coefficient can be nonzero only when the total Weyl index is the
negative of that of $PQ$; we therefore need only consider
$\theta_{1234}=1$. If all four matrices have the same type, every
$t_{ij}$ is one and \eqref{eq:H}--\eqref{eq:b} give $r_{1234}=0$.
The other cases use the following matrices and the indicated transpose:
\begin{align*}
 M_{31}&=\begin{pmatrix}
 1&1&1&v\\1&1&1&w\\1&1&1&(vw)^{-1}\\u&z&(uz)^{-1}&1
 \end{pmatrix},&
 M_{\mathrm{opp}}&=\begin{pmatrix}
 1&1&v&w\\1&1&v^{-1}&w^{-1}\\u&u^{-1}&1&1\\z&z^{-1}&1&1
 \end{pmatrix},\\[4pt]
 M_{\mathrm{both}}&=\begin{pmatrix}
 1&1&v&v^{-1}\\1&1&v^{-1}&v\\u&u^{-1}&1&1\\u^{-1}&u&1&1
 \end{pmatrix},&
 M_{\mathrm{split}}&=\begin{pmatrix}
 1&1&v&v^{-1}\\1&1&w&w^{-1}\\u&u^{-1}&1&1\\z&z^{-1}&1&1
 \end{pmatrix}.
\end{align*}
Here $u,z$ are nontrivial roots of one of the two prime orders, and
$v,w$ of the other; the prime orders may be interchanged. For $M_{31}$
and $M_{\mathrm{split}}$ (including its transpose), we require
$uz\ne1$ and $vw\ne1$. For $M_{\mathrm{opp}}$, we require only
$(uz,vw)\ne(1,1)$; the excluded case is $M_{\mathrm{both}}$.
The triple characters and the resulting quartic coefficients are
\begin{equation}\label{eq:zero-table}
\renewcommand{\arraystretch}{1.25}
\begin{array}{c|cccc|c}
 (t_{ij})&\theta_{123}&\theta_{124}&\theta_{134}&\theta_{234}&r_{1234}\\\hline
 M_{31}&1&uzvw&(zw)^{-1}&(uv)^{-1}&0\\
 M_{\mathrm{opp}}&1&1&uzvw&(uzvw)^{-1}&0\\
 M_{\mathrm{both}}&1&1&1&1&0\\
 M_{\mathrm{split}},\ M_{\mathrm{split}}^{\mathsf T}
   &vw&(vw)^{-1}&uz&(uz)^{-1}&0
\end{array}
\end{equation}
To obtain these identities, first set $b_{ij}=0$ for same-type pairs and
$b_{ij}=(2-t_{ij}-t_{ji})/(1-t_{ij}t_{ji})$ for mixed pairs,
then apply \eqref{eq:H}--\eqref{eq:b} for triples and for the full set.

These cases exhaust all possibilities. With three matrices of one type,
their two coordinate sums must both be zero, giving $M_{31}$ after
placing the repeated type first. With two matrices of each type,
the conditions are
\[
 (a_1+a_2)(c_1+c_2)=0,\qquad
 (b_1+b_2)(d_1+d_2)=0.
\]
If one pair is opposite in both coordinates, place it first and use
$M_{\mathrm{opp}}$, or $M_{\mathrm{both}}$ if both pairs are opposite.
Otherwise the zero sums occur in different pairs in the two coordinates,
giving $M_{\mathrm{split}}$ or its transpose.
A product of a $p$-th root and a $q$-th root is one only when both roots
are one. Thus the classification above identifies exactly the subsets $I$
for which $\theta_I=1$, and every division in the recurrence is
valid after specialization. This proves the proposition.
\end{proof}

\subsection{The fifth order condition}

In this section we show that the fifth order condition fails for $n$ a product of two distinct primes.

For a commuting basis pair $P,Q$, let $\Omega_{P,Q}$ denote the fifth homogeneous term of
$\mathcal R_{P,Q}$. Equation~\eqref{eq:Rquartic} and
Proposition~\ref{prop:fourzero} give
\begin{equation}\label{eq:leading}
 \mathcal R_{P,Q}(X)=\Omega_{P,Q}(X)+O(\|X\|_2^6).
\end{equation}

\begin{lemma}\label{lem:fifth}
Set $u=e^{2\pi i/p}$, $v=e^{2\pi i/q}$, and choose
\begin{align*}
 W_1&=P_p\otimes Q_q,&W_2&=P_p\otimes Q_q^{-2},&
 W_3&=P_p^{-2}\otimes Q_q,\\
 W_4&=Q_p\otimes P_q,&W_5&=W_4.
\end{align*}
For $P_0=I_p\otimes P_q^{-2}$ and $Q_0=Q_p^{-2}\otimes I_q$,
the coefficient of $z_1z_2z_3z_4z_5$ in
$\Omega_{P_0,Q_0}(\sum_i z_iW_i)$ is
\begin{equation}\label{eq:fifth}
 \frac{4(u-1)^6(u+1)^2(v-1)^6(v+1)^2}
 {(u^2-v^4)(u^4-v^2)(u^2v^2-1)}.
\end{equation}
\end{lemma}
\begin{proof}
The first three Weyl indices sum to zero. The total index is therefore
$((2,0),(0,2))$, and $P_0,Q_0$ commute. The cross-character matrix is
\begin{equation}\label{eq:fifth-matrix}
 (t_{ij})=\begin{pmatrix}
 1&1&1&v&v\\
 1&1&1&v^{-2}&v^{-2}\\
 1&1&1&v&v\\
 u&u&u^{-2}&1&1\\
 u&u&u^{-2}&1&1
 \end{pmatrix}.
\end{equation}
Among subsets $I$ with at least two elements, the associated Weyl
matrix $E_I$ belongs to $\T$ precisely when $I$ is a subset of
$\{1,2,3\}$, the pair $\{4,5\}$, one of the sets
$\{1,2,3,4\}$ and $\{1,2,3,5\}$, or the full set.
For every other such subset, $\chi(x_I,y_I)\ne1$, so $E_I\in\Sspace$. Indeed, the sum of
one or two of the first three indices has a nonzero $P_p$-coordinate
and a nonzero $Q_q$-coordinate. One or two copies of $W_4$ have the same
property in the complementary coordinates. This uses that $p,q$ are odd.

Apply \eqref{eq:H}--\eqref{eq:b} to \eqref{eq:fifth-matrix}. For clarity,
write $a_i=t_{4i}=t_{5i}$ and $c_i=t_{i4}=t_{i5}$, so that
\[
 (a_1,a_2,a_3)=(u,u,u^{-2}),\qquad
 (c_1,c_2,c_3)=(v,v^{-2},v).
\]
The first mixed coefficients are
\begin{equation}\label{eq:fifth-intermediate}
 b_{i4}=b_{i5}=\frac{a_i+c_i-2}{a_ic_i-1},\qquad
 b_{i45}=\frac{2(a_i-1)^2(c_i-1)^2}
 {(a_ic_i-1)^2(a_ic_i+1)},\quad 1\le i\le3.
\end{equation}
In particular, the factor $2$ accounts for the two independent labels of
the repeated matrix $W_4=W_5$. Put
\[
 \Delta=(u-v^2)^2(u+v^2)(u^2-v)^2(u^2+v)(uv-1)^2(uv+1).
\]
Substitution and collection over this common denominator give
\begin{align*}
 \Delta r_{12345}={}&4(u-1)^6(u+1)^2(v-1)^6(v+1)^2\\
 &\hspace{16mm}\cdot(u-v^2)(u^2-v)(uv-1).
\end{align*}
This is an identity of Laurent polynomials after clearing denominators
in the recurrence. Cancelling the three displayed factors gives
\eqref{eq:fifth}. The formal labels $z_i$ are independent, so the
calculation remains valid when Weyl matrices coincide, as happens for
$W_4,W_5$, and also when one prime is $3$.
\end{proof}

\begin{proposition}\label{prop:odd}
If $n=pq$ with $p,q$ distinct odd primes, there exists a nonzero
$a\in\A_n$ such that $a/\|a\|_2$ is not a direction of convergence.
\end{proposition}
\begin{proof}
Every numerator factor in \eqref{eq:fifth} is nonzero, since $p,q$ are
odd. The denominator is nonzero as well. For example, $u^2=v^4$ would
be a root of unity of both $p$-power and $q$-power order, forcing
$u^2=1$. The other two factors are treated in the same way.
Thus $\Omega_{P_0,Q_0}$ is a nonzero polynomial on $\T$.
Since $\A$ is a real form of $\T$, as in the proof of
Proposition~\ref{prop:two}, some $a\in\A$ satisfies
$\Omega_{P_0,Q_0}(a)\ne0$.

Suppose that $a/\|a\|_2$ is a direction of convergence. Write
$V_k=U_k-I=X_k+N_k$ in $\T\oplus\Sspace$, and put
$\varepsilon_k=\|V_k\|_2$. The implicit function theorem gives
$N_k=Y(X_k)$ for large $k$. By \eqref{eq:leading},
\[
 0=\Omega_{P_0,Q_0}(X_k)+O(\|X_k\|_2^6).
\]
Since $X_k/\varepsilon_k\to a/\|a\|_2$ and
$\|X_k\|_2\le\varepsilon_k$, division by $\varepsilon_k^5$ gives
$\Omega_{P_0,Q_0}(a/\|a\|_2)=0$, a contradiction. The Chinese remainder
identification carries the conclusion back to $F_n$.
\end{proof}

\begin{remark}
The oddness assumption is used both in the subset classification in
Lemma~\ref{lem:fifth} and in the factors $u+1,v+1$ of the obstruction.
When one prime is $2$, the equality $\theta_I=1$ holds for some
additional subsets $I$, so the rational calculation cannot simply
be specialized to that case.
\end{remark}

\section{The main result}\label{sec:consequences}

In this section we prove our main result, then present some consequences.

\begin{proof}[Proof of Theorem~\ref{thm:main}]
Under the hypothesis, $n$ has at least three distinct prime divisors,
has exactly two with a repeated factor, or is the product of two distinct
odd primes. Apply Propositions~\ref{prop:three}, \ref{prop:two},
and~\ref{prop:odd}, respectively.
\end{proof}

\begin{corollary}\label{cor:dimension}
Assume that $n\ge2$ is neither a prime power nor twice a prime.
Let $\Omega\subset\R^m$ be open, $0\in\Omega$, and let
$f:\Omega\to\Hset(n)$ be differentiable at $0$, with
$f(0)=\sqrt n F_n$. Then
\[
 \rank_{\R}Df(0)\le d_n-1.
\]
If the matrices $f(x)$ are dephased, then
$\rank_{\R}Df(0)\le d_n'-1=d_n-2n$.
\end{corollary}
\begin{proof}
Set $g(x)=n^{-1/2}f(x)F_n^*$. Then $g(0)=I_n$ and
$\operatorname{im}Dg(0)\subset\A_n$. If $\rank Dg(0)=d_n$, this image
is all of $\A_n$. Choose $a$ as in Theorem~\ref{thm:main} and
$h\in\R^m$ such that $Dg(0)h=a$. For $\lambda\downarrow0$,
\[
 g(\lambda h)=I_n+\lambda a+o(\lambda),
\]
which realizes the excluded direction $a/\|a\|_2$, a contradiction.
For a dephased family of rank $k$, adjoining $n$ row phases and $n-1$
column phases gives a family of rank $k+2n-1$. The added tangent space meets
the dephased tangent space trivially, as follows by inspecting the first
row and column. Thus $k+2n-1\le d_n-1$.
\end{proof}

For $n=pq$ with $p,q$ distinct odd primes, \eqref{eq:defect} gives
\[
 d_{pq}=(2p-1)(2q-1),\qquad d_{pq}'=2(p-1)(q-1).
\]
Thus the derivative of a dephased family through $\sqrt{pq}F_{pq}$ has
rank at most $2(p-1)(q-1)-1$. In particular, at $\sqrt{15}F_{15}$ its
rank is at most $15$, whereas the dephased defect is $16$.

\end{document}